\documentclass[11pt]{article}
\usepackage{amsmath,amsfonts,amssymb,amsthm}
\begin{document}

\newcommand{\1}{{{\bf 1}}}
\newcommand{\id}{{\rm id}}
\newcommand{\Hom}{{\rm Hom}\,}
\newcommand{\End}{{\rm End}\,}
\newcommand{\Res}{{\rm Res}\,}
\newcommand{\Image}{{\rm Im}\,}
\newcommand{\Ind}{{\rm Ind}\,}
\newcommand{\Aut}{{\rm Aut}\,}
\newcommand{\Ann}{{\rm Ann}}
\newcommand{\Ker}{{\rm Ker}\,}
\newcommand{\gr}{{\rm gr}}
\newcommand{\Der}{{\rm Der}\,}

\newcommand{\Z}{\mathbb{Z}}
\newcommand{\Q}{\mathbb{Q}}
\newcommand{\C}{\mathbb{C}}
\newcommand{\N}{\mathbb{N}}
\newcommand{\g}{\mathfrak{g}}
\newcommand{\h}{\mathfrak{h}}
\newcommand{\wt}{{\rm wt}\;}
\newcommand{\A}{\mathcal{A}}
\newcommand{\D}{\mathcal{D}}
\newcommand{\Lie}{\mathcal{L}}

\def \<{\langle}
\def \>{\rangle}
\def \be{\begin{equation}\label}
\def \ee{\end{equation}}
\def \bex{\begin{exa}\label}
\def \eex{\end{exa}}
\def \bl{\begin{lem}\label}
\def \el{\end{lem}}
\def \bt{\begin{thm}\label}
\def \et{\end{thm}}
\def \bp{\begin{prop}\label}
\def \ep{\end{prop}}
\def \br{\begin{rem}\label}
\def \er{\end{rem}}
\def \bc{\begin{coro}\label}
\def \ec{\end{coro}}
\def \bd{\begin{de}\label}
\def \ed{\end{de}}

\newtheorem{thm}{Theorem}[section]
\newtheorem{prop}[thm]{Proposition}
\newtheorem{coro}[thm]{Corollary}
\newtheorem{conj}[thm]{Conjecture}
\newtheorem{exa}[thm]{Example}
\newtheorem{lem}[thm]{Lemma}
\newtheorem{rem}[thm]{Remark}
\newtheorem{de}[thm]{Definition}
\newtheorem{hy}[thm]{Hypothesis}
\makeatletter \@addtoreset{equation}{section}
\def\theequation{\thesection.\arabic{equation}}
\makeatother \makeatletter

\begin{flushright}
\today
\\
\end{flushright}

\begin{Large}
\begin{center}
\textbf{Twisted modules for vertex algebras associated with vertex algebroids}
\end{center}
\end{Large}
\begin{center}{Haisheng Li\footnote{Partially supported by an NSA grant}\\
Department of Mathematical Sciences, Rutgers University, Camden, NJ 08102\\
and\\
Department of Mathematics, Harbin Normal University, Harbin, China\\
\ \
\\
Gaywalee Yamskulna\\
Department of Mathematics, Illinois State University, Normal, IL
61790\\
and\\
Institute of Science, Walailak University, Nakhon Si Thammarat,
Thailand}
\end{center}

\begin{abstract}
We continue with \cite{liy} to construct and classify graded simple twisted
modules for the $\N$-graded vertex algebras constructed by
Gorbounov, Malikov and Schechtman from vertex algebroids.
Meanwhile  we determine the
full automorphism groups of those $\N$-graded vertex algebras
in terms of the automorphism groups of the corresponding vertex algebroids.
\end{abstract}

\section{Introduction}
For most of the important examples of vertex operator algebras
$V=\coprod_{n\in \Z}V_{(n)}$ graded by the $L(0)$-weight (see [FLM,
FHL]), the $\Z$-grading satisfies the condition that $V_{(n)}=0$ for
$n<0$ and $V_{(0)}=\C {\bf 1}$ where ${\bf 1}$ is the vacuum vector.
For a vertex operator algebra $V$ with this special property, the
homogeneous subspace $V_{(1)}$ has a natural Lie algebra structure
with $[u,v]=u_{0}v$ for $u,v\in V_{(1)}$ and the product $u_{1}v$
$(\in V_{(0)})$ defines a symmetric invariant bilinear form on
$V_{(1)}$.

In a series of study on
Gerbs of chiral differential operators in [GMS] and on
chiral de Rham complex in [MSV, MS1,2], Malikov and his coauthors investigated
$\N$-graded vertex algebras $V=\coprod_{n\in \N}V_{(n)}$ with $V_{(0)}$
{\em not} necessarily $1$-dimensional.  
In this case, the bilinear operations $(u,v)\mapsto u_{i}v$ for $i\ge 0$
are closed on $V_{(0)}\oplus V_{(1)}$:
$$u_{i}v\in V_{(0)}\oplus V_{(1)} \;\;\;\mbox{ for }u,v\in
V_{(0)}\oplus V_{(1)},\; i\ge 0.$$ The skew symmetry and the Jacobi
identity for the vertex algebra $V$ give rise to several compatibility
relations.  Such algebraic structures on $V_{(0)}\oplus V_{(1)}$ are
summarized in the notion of what was called a $1$-truncated conformal
algebra. Furthermore, the subspace $V_{(0)}$ equipped with the product
$(a,b)\mapsto a_{-1}b$ is a commutative associative algebra with the
vacuum vector ${\bf 1}$ as the identity and $V_{(0)}$ as a
nonassociative algebra acts on $V_{(1)}$ by $a\cdot u=a_{-1}u$ for
$a\in V_{(0)}$, $u\in V_{(1)}$. All these structures on $V_{(0)}\oplus
V_{(1)}$ are further summarized in the notion of what was called a
vertex $A$-algebroid, where $A$ is a (unital) commutative associative
algebra.  On the other hand, in [GMS], among other important results,
Gorbounov, Malikov and Schechtman constructed an $\N$-graded vertex
algebra $V=\coprod_{n\in \N}V_{(n)}$ {}from any vertex $A$-algebroid,
such that $V_{(0)}=A$ and the vertex $A$-algebroid $V_{(1)}$ is
isomorphic to the given one.  All the constructed $\N$-graded vertex
algebras are generated by $V_{(0)}\oplus V_{(1)}$ with a spanning
property of PBW type.  As it was demonstrated in [GMS], such
$\N$-graded vertex algebras are natural and important to study. For
example, the vertex (operator) algebra associated with a $\beta\gamma$
system, which plays a central role in free field realization of affine
Lie algebras (see [W, FF1-3, FB]) is such an $\N$-graded vertex
algebra.  The vertex (operator) algebras constructed from toroidal Lie
algebras are also of this type (see [BBS, BDT]).

In \cite{liy}, we revisited those $\N$-graded vertex algebras
and we classified all the $\N$-graded simple modules
in terms of simple modules for certain Lie algebroids.
In the theory of vertex algebras, in addition to the notion of module
we have the notion of twisted module and twisted modules play 
a very important role, especially in the study of the so-called
orbifold theory. Certainly, twisted modules also play an important role 
in other studies.
In this paper, we continue to study the twisted modules 
for the $\N$-graded vertex algebras associated with vertex algebroids.

Let $B$ be a vertex $A$-algebroid and let $V_{B}$ be the associated
$\N$-graded vertex algebra.  In this paper, we define a notion of
automorphism of the vertex $A$-algebroid $B$ and we prove that any
automorphism of the vertex $A$-algebroid $B$ can be extended uniquely
to an automorphism of the $\N$-graded vertex algebra $V_{B}$ and that
the full automorphism group of the $\N$-graded vertex algebra $V_{B}$
is naturally isomorphic to the full automorphism group of the vertex
$A$-algebroid $B$.  Let $g$ be an automorphism of the vertex
$A$-algebroid $B$ of order $T$ (finite).  Then the $g$-fixed point
$A^{0}$ is a subalgebra of $A$ and the $g$-fixed point $B^{0}$ is a
vertex $A^{0}$-algebroid.  Furthermore, $B^{0}/A^{0}\partial A^{0}$ is
a Lie $A^{0}$-algebroid. It is proved that the category of
$\frac{1}{T}\N$-graded simple $g$-twisted $V_{B}$-modules is
equivalent to a subcategory of simple modules for the Lie
$A^{0}$-algebroid $B^{0}/A^{0}\partial A^{0}$.

This paper is organized as follows:
In Section 2, we review the construction of vertex algebras 
associated with vertex algebroids and we identify their automorphism groups
with the automorphism groups of the vertex algebroids.
In Section 3, we classify graded simple twisted modules.

\section{Preliminaries}

We recall the notions of $1$-truncated conformal algebra, vertex
algebroid and Lie algebroid, and we review the construction
of the $\N$-graded vertex algebra $V_{B}$ associated with a vertex $A$-algebroid $B$.
We also define notions of (endomorphism) automorphism of 
a $1$-truncated conformal algebra and of a vertex
$A$-algebroid $B$. We then identify the group of 
grading-preserving automorphisms of
$V_B$ with the group of automorphisms of the vertex $A$-algebroid $B$.

First, we recall from [GMS] (cf. [Br1-2]) the notions of $1$-truncated conformal algebra,
vertex algebroid and Lie algebroid.

\begin{de}
{\em A {\em $1$-truncated conformal algebra} is a
graded vector space $C=C_0\oplus C_1$, equipped with a linear map
$\partial: C_0\rightarrow C_1$ and bilinear operations
$(u,v)\mapsto u_iv$ for $i=0, 1$ of degree $-i-1$ on $C$ such that
the following axioms hold:
\begin{enumerate}
\item (Derivation) for $a\in C_0$, $u\in C_1$,
\begin{equation}
(\partial a)_0=0;\ \ (\partial a)_1=-a_0;\ \ \partial(u_0a)=u_0\partial a
\end{equation}
\item (Commutativity) for $a\in C_0$, $u,v\in C_1$,
\begin{equation}
u_0a=-a_0u;\ \ u_0 v=-v_0u+\partial (v_1 u);\ \ u_1v=v_1u
\end{equation}
\item (Associativity) for $\alpha , \beta ,\gamma \in C_0\oplus C_1$,
\begin{equation}
\alpha _0 \beta _i \gamma=\beta _i \alpha _0 \gamma +(\alpha _0\beta)_i\gamma.
\end{equation}
\end{enumerate}}
\end{de}


\br{rtca-deform}
{\em Let $C=C_{0}\oplus C_{1}$ be a $1$-truncated conformal
algebra and let $\ell$ be any nonzero complex number.
Set $C[\ell]=C_{0}\oplus C_{1}$ as a vector space.
We retain all the structures on $C$ except that 
we change the bilinear operation $C_{1}\times C_{1}\rightarrow C_{0}:
u\times v\mapsto u_{1}v$
by multiplying $1/\ell$ and change the linear operator $\partial$ by
multiplying $\ell$.
Then one can show that $C[\ell]$ is a $1$-truncated conformal algebra.}
\er

\begin{de}\label{def-vad} 
{\em Let $A$ be a unital commutative associative algebra over $\C$.
A {\em vertex $A$-algebroid} is a $\C$-vector space $\Gamma$
equipped with
\begin{enumerate}
\item a $\C$-bilinear map $$A \times \Gamma \rightarrow \Gamma ; \ \  (a,v)\mapsto a*v$$
such that $1*v=v$ for $v\in \Gamma$.
\item a structure of a Leibniz 
$\C$-algebra $[\cdot,\cdot ]:\Gamma\otimes_{\C}\Gamma\rightarrow \Gamma$.
\item a homomorphism of Leibniz $\C$-algebras $\pi : \Gamma \rightarrow Der(A)$.
\item a symmetric $\C$-bilinear pairing $\<\cdot,\cdot\>:\Gamma \otimes_{\C}\Gamma \rightarrow A$.
\item a $\C$-linear map $\partial: A\rightarrow \Gamma$ such that $\pi\circ \partial =0$.
\end{enumerate}
All the following conditions are assumed to hold:
\begin{eqnarray*}
a*(a'*v)-(aa')*v&=&\pi(v)(a)*\partial (a')+\pi (v)(a')*\partial (a)\\
{[u,a*v]}&=&\pi (u)(a)*v+a*[u,v]\\
{[u,v]+[v,u]}&=&\partial (\<  u,v\>)\\
{\pi(a*v)}&=&a\pi(v)\\
{\<a*u,v\>}&=&a\<u,v\>-\pi(u)(\pi(v)(a))\\
{\pi(v)(\<v_1,v_2\>)} &=&\<[v,v_1],v_2\>+\<v_1,[v,v_2]\>\\
{\partial (aa')}&=&a*\partial (a')+a'*\partial (a)\\
{[v,\partial (a)]}&=&\partial (\pi(v)(a))\\
{\< v,\partial (a)\>}&=&\pi(v)(a)
\end{eqnarray*}
for $a,a'\in A$, $u,v,v_1,v_2\in \Gamma.$}
\end{de}


The following was proved in \cite{liy}:

\begin{prop}\label{pLiY} 
Let $A$ be a unital commutative associative algebra and let $B$ be a module for $A$ 
as a nonassociative algebra.
Then a vertex $A$-algebroid structure on $B$ is equivalent to 
a $1$-truncated conformal algebra structure on $C=A\oplus B$ with
\begin{eqnarray}
& &a_ia'=0,\\
& &u_0v=[u,v],\ \ u_1v=\< u,v\>,\\
& &u_0a=\pi(u)(a),\ \ a_0u=-u_0a=-\pi(u)(a)
\end{eqnarray}
for $a,a'\in A$, $u,v\in B$, $i=0,1$, such that
\begin{eqnarray}
& &a(a'u)-(aa')u=(u_0a)\partial a'+(u_0a')\partial a,\\
& &u_0(av)-a(u_0v)=(u_0a)v,\\
& &u_0(aa')=a(u_0a')+(u_0a)a',\\
& &a_0(a'v)=a'(a_0v),\\
& & (au)_1v=a(u_1v)-u_0v_0a,\\
& &\partial (aa')=a\partial (a')+a' \partial (a).
\end{eqnarray}
\end{prop}

\begin{de} 
{\em Let $A$ be a unital commutative associative algebra. A {\em
Lie $A$-algebroid} is a Lie algebra $\g$ equipped with an
$A$-module structure and a module action on $A$ by derivation such
that
\begin{eqnarray}
[u,av]&=&a[u,v]+(ua)v,\\
a(ub)&=&(au)b\ \ \text{ for }u,v\in \g,\; a,b\in A.
\end{eqnarray}

A {\em module} for a Lie $A$-algebroid $\g$ is a vector space
$W$ equipped with a $\g$-module structure and an $A$-module
structure such that
\begin{eqnarray}
& &u(aw)-a(uw)=(ua)w,\\
& &a(uw)=(au)w \ \ \ \text{ for }a\in A,\; u\in \g,\; w\in W.
\end{eqnarray}}
\end{de}

The following result was due to [Br2]:

\begin{lem}\label{br2} 
Let $A$ be a unital commutative associative algebra
(over $\C$) and let $B$ be a vertex $A$-algebroid. Then
$B/A\partial A$ is naturally a Lie $A$-algebroid.
\end{lem}

Next, we recall the construction of vertex algebras associated
with vertex algebroids, following the exposition of \cite{liy}.

First, starting with a $1$-truncated conformal algebra $C=A\oplus B$ 
we construct a Lie algebra.
Set
\begin{eqnarray}
L(A\oplus B)=(A\oplus B)\otimes \C[t,t^{-1}].
\end{eqnarray}
In the obvious way we define the subpaces $L(A)$ and $L(B)$.
Set
\begin{eqnarray*}
\hat\partial =\partial\otimes 1+1\otimes d/dt: \;\;
L(A)\rightarrow L(A\oplus B).
\end{eqnarray*}
We define
\begin{eqnarray*}
\deg (a\otimes t^{n})&=&-n-1\;\;\;\mbox{ for }a\in A,\;n\in \Z,\\
\deg (b\otimes t^{n})&=&-n\;\;\;\mbox{ for }b\in B,\;n\in \Z,
\end{eqnarray*}
making $L(A\oplus B)$ a $\Z$-graded vector space. The linear
map $\hat\partial$ is homogeneous of degree $1$.
Set 
\begin{eqnarray}
\Lie=L(A\oplus B)/\hat\partial L(A).
\end{eqnarray}

Define a bilinear product $[\cdot,\cdot]$ on $L(A\oplus B)$ such that
for $a,a'\in A,\; b,b'\in B,\;m,n\in \Z$,
\begin{eqnarray}
& &[a\otimes t^{m}, a'\otimes t^{n}]=0,\label{edef-bracket-1}\\
& &{[a\otimes t^{m}, b\otimes t^{n}]}=a_{0}b\otimes t^{m+n},\label{edef-bracket-2}\\
& &{[b\otimes t^{n},a\otimes t^{m}]}=b_{0}a\otimes
t^{m+n},\label{edef-bracket-22}\\
& &{[b\otimes t^{m}, b'\otimes t^{n}]}=b_{0}b'\otimes
t^{m+n}+m(b_{1}b')\otimes t^{m+n-1}.\label{edef-bracket-3}
\end{eqnarray}
The following result was established in \cite{liy}:

\bp{pto-quotient}
Let $C=A\oplus B$ be a $1$-truncated
conformal algebra. The subspace $\hat{\partial}L(A)$ of the nonassociative
algebra $(L(A\oplus B),[\cdot,\cdot])$ is a two-sided ideal. 
Furthermore, the quotient nonassociative algebra $\Lie$ is a $\Z$-graded 
Lie algebra.
\ep 

Let $\rho$ be the projection map from $L(A\oplus B)$ to
$\Lie$. For $u\in A\oplus B,\; n\in \Z$, we set
$$u(n)=\rho(u\otimes t^n)=u\otimes t^n+\hat\partial L(A)\in \Lie.$$ 
We have graded Lie subalgebras
\begin{eqnarray*}
\Lie^{\ge 0}&=&\rho ((A\oplus B)\otimes \C[t]),\\
\Lie^{< 0}&=&\rho ((A\oplus B)\otimes t^{-1}\C[t^{-1}])
\end{eqnarray*}
and we have $\Lie=\Lie^{\ge 0}\oplus \Lie^{<0}$ as a
vector space.

Considering $\C$ as a trivial $\Lie^{\ge 0}$-module we
form the induced module
\begin{equation*}
V_{\Lie}=U(\Lie)\otimes _{U(\Lie^{\ge 0})}\C.
\end{equation*}
We assign $\deg \C=0$, making $V_{\Lie}$ naturally an
$\N$-graded $\Lie$-module:
\begin{eqnarray}\label{egrading-valie}
V_{\Lie}=\coprod_{n\in\N}(V_{\Lie})_{(n)}.
\end{eqnarray}
Throughout this paper, $\N$ denotes the set of nonnegative integers. 
Set
$${\bf 1}=1\otimes 1\in V_{\Lie}.$$
By the P-B-W theorem, we have $V_{\Lie}=U(\Lie^{< 0})=S(\Lie^{< 0})$.
In view of this, we can and we do consider $A\oplus B$ as a subspace:
$$A\oplus B\rightarrow V_{\Lie};\;\; u\mapsto u(-1){\bf 1}.$$ 

The following was proved in \cite{liy} (cf. \cite{dlm-vpa}):

\bt{tlie-vertex-algebra}
There exists a unique vertex algebra structure on $V_{\Lie}$
with ${\bf 1}$ as the vacuum vector and with
$Y(u,x)=\sum_{n\in\Z}u(n)x^{-n-1}$ for $u\in A\oplus B$. Moreover,
the vertex algebra $V_{\Lie}$ is naturally an $\N$-graded vertex algebra and 
is generated by the subspace $A\oplus B$ with $A$ of degree $0$ and $B$ of degree $1$.
\et

\br{rbasis}
{\em For $n\in \Z$, set
$$A(n)=\{a(n)\; |\; a\in A\}, \ \ B(n)=\{b(n)\; |\;  b\in B\}\subset \Lie,$$ 
and we set
$$B(-)= \coprod_{n=1}^{\infty}B(-n)\subset \Lie.$$ 
Both $A(-1)$ and $B(-)$ are Lie
subalgebras of $\Lie ^{<0}$ and we have ${\Lie}^{<0}=A(-1)\oplus B(-)$
as a vector space.
Then
\begin{eqnarray}
V_{\Lie}=U({\Lie}^{<0})=S({\Lie}^{<0})
=S(A(-1)\oplus B(-))=S(B(-))\otimes S(A(-1)).
\end{eqnarray}
Consequently, $(V_{\Lie})_{(n)}=S(B(-))_{(n)}\otimes S(A(-1))$ for $n\in \N$.
In particular, $(V_{\Lie})_{(0)}= S(A(-1))$. }
\er

Now, we assume that $A$ is a unital
commutative associative algebra with the identity $e$ and $B$ is a
vertex $A$-algebroid. In particular, $C=A\oplus B$ is a $1$-truncated
conformal algebra. We set 
\begin{eqnarray*} 
E&=&{\rm span}\{e-{\bf 1},\ \
a(-1)a'-aa',\ \ a(-1)b-ab\; |\; a,a'\in A, \; b\in B\}\subset V_{\Lie},\\
I_B&=&U(\Lie)\C [\D]E.
\end{eqnarray*}
It was proved in \cite{liy} that the $\Lie$-submodule $I_B$ of $V_{\Lie}$ 
is a two-sided graded ideal of the $\N$-graded vertex algebra $V_{\Lie}$.
The $\N$-graded vertex algebra $V_{B}$ associated with the vertex
$A$-algebroid $B$ is defined to be the quotient vertex algebra
\begin{eqnarray}
V_B=V_{\Lie}/I_B.
\end{eqnarray}
We have (see [GMS], \cite{liy}):

\begin{prop}\label{pgms-ly}
Let $A$ be a unital commutative associative algebra with the identity $e$ and $B$ a
vertex $A$-algebroid. Then $V_B$ is an $\N$-graded vertex algebra such that
$(V_B)_{(0)}=A$, $(V_B)_{(1)}=B$ and for $n\geq 1$,
\begin{eqnarray*}
& &(V_B)_{(n)}\\
&=&{\rm span}\{b_1(-n_1)\cdots b_k(-n_k){\bf 1}\;|\;b_i\in B, n_1\geq n_2\geq
\cdots\geq n_k\geq 1, n_1+\cdots +n_k=n\}.
\end{eqnarray*}
In particular, $V_B$ is generated by the subspace $A\oplus B$.
\end{prop}

Next, we discuss homomorphisms and automorphisms for
$1$-truncated conformal algebras, vertex $A$-algebroids and 
for the $\N$-graded vertex algebras $V_{B}$.

\bd{d-end-ca}
{\em Let $C=A\oplus B$ and $C'=A'\oplus B'$ be $1$-truncated conformal
algebras. A {\em homomorphism} from $C$ to $C'$ is a linear map
$f: C\rightarrow C'$ such that $f(A)\subset A',\; f(B)\subset B'$,
$f\partial=\partial f$,  and such that
$$f(u_{i}v)=f(u)_{i}f(v)$$
for $u,v\in C,\; i=0,1$.}
\ed

\begin{lem}\label{lend-lie} 
Let $f$ be an endomorphism of a $1$-truncated conformal
algebra $C=A\oplus B$.
Then the linear endomorphism of $L(A\oplus B)$ defined by
\begin{eqnarray}
\hat{f}(u\otimes t^{n})=f(u)\otimes t^{n}
\end{eqnarray}
for $u\in A\oplus B,\; n\in \Z$ gives rise to
an endomorphism of $\Lie$, which we denote by $\hat{f}$ again.
Furthermore, $\hat{f}$ preserves the $\Z$-grading of $\Lie$.
\end{lem}

\begin{proof} Using the property that $f\partial=\partial f$, 
we have $\hat{f}\hat{\partial}=\hat{\partial}\hat{f}$. For $u,v\in C=A\oplus B$,
as $f(u_{i}v)=f(u)_{i}f(v)$ for $i=0,1$, from 
(\ref{edef-bracket-1})-(\ref{edef-bracket-3}) 
we have
$$\hat{f}\left([u\otimes t^{m},v\otimes t^{n}]\right)
=[f(u)\otimes t^{m},f(v)\otimes t^{n}]
=[\hat{f}(u\otimes t^{m}),\hat{f}(v\otimes t^{n})].$$
Thus $\hat{f}$ gives rise to an endomorphism of the Lie algebra $\Lie$.
It is clear that $\hat{f}$ preserves the $\Z$-grading.
\end{proof}

\begin{de} {\em Let $A$ and $A'$ be unital commutative associative algebras and
let $B$ be a vertex $A$-algebroid, $B'$ a vertex $A'$-algebroid. 
A {\em vertex algebroid homomorphism} from $B$ to $B'$
is a linear map $f: A\oplus B\rightarrow A'\oplus B'$ 
such that $f(A)\subset A',\; f(B)\subset B'$ and such that
\begin{enumerate}
\item $f|_A$ is an associative algebra homomorphism. 
\item $f|_B$
is a Leibniz algebra homomorphism. 
\item $f(ab)=f(a)f(b)$ for
$a\in A,\; b\in B$. 
\item $\<f(u),f(v)\>=f(\<u,v\>)$ for $u,v\in B$. 
\item $f\circ \partial =\partial \circ f$.
\item $f(b_0a) =f(b)_0f(a)$ for $a\in A,\;b\in B$.
\end{enumerate}
An {\em automorphism of a vertex $A$-algebroid $B$} is
a bijective vertex algebroid endomorphism of the vertex $A$-algebroid $B$.}
\end{de}

Let $(V,Y, {\bf 1})$ be a vertex algebra. An {\em endomorphism} of $V$
is a linear map $g: V\rightarrow V$ such that
\begin{eqnarray}
g({\bf 1})&=&{\bf 1},\label{au1}\\
g(Y(u,x)v)&=&Y(g(u),x)g(v)\label{au2}
\end{eqnarray}
for $u,v\in V$. An {\em automorphism} of $V$
is a bijective endomorphism of $V$. 
The group of automorphisms of $V$ is denoted by $\Aut(V)$.
If $V=\coprod_{m\in\Z}V_{(m)}$ is a $\Z$ (or $\N$)-graded vertex algebra,
we denote by $\Aut^{o}(V)$ the group of 
grading-preserving automorphisms of $V$.

\begin{lem}\label{lav} 
Let $B$ be a vertex $A$-algebroid and let $g$ be a grading-preserving automorphism 
of the vertex algebra $V_B$.
Then $g$ restricted to $A\oplus B$ is an automorphism of the vertex $A$-algebroid $B$.
\end{lem}

\begin{proof} As $(V_{B})_{(0)}=A$ and $(V_{B})_{(1)}=B$,
$g$ is a linear bijection on $A\oplus B$ that preserves the
subspaces $A$ and $B$. For $a,a'\in A, \;b,b'\in B$, we have
\begin{eqnarray*}
& &g(aa')=g(a(-1)a')=g(a)_{-1}g(a')=g(a)g(a'),\\
& &g(ab)=g(a(-1)b)=g(a)_{-1}g(b)=g(a)g(b),\\
& &g([b,b'])=g(b_0b')=g(b)_0g(b')=[g(b),g(b')],\\
& &g(\<b,b'\>)=g(b_1b')=g(b)_1g(b')=\<g(b),g(b')\>,\\
& &g(b_0a)=g(b)_0g(a),\\ 
& &g(\partial(a))=g(a(-2){\bf 1})=g(a)_{-2}{\bf 1}=\partial(g(a)).
\end{eqnarray*}
Thus $g$ is an automorphism of vertex $A$-algebroid $B$.
\end{proof}

On the other hand, we are going to prove that any
automorphism of a vertex $A$-algebroid $B$ extends canonically
to an automorphism of the $\N$-graded vertex algebra $V_{B}$.
First we have:

\bl{lva} 
Let $C=A\oplus B$ be a $1$-truncated conformal algebra and let
$g$ be an endomorphism of $C$.
Then $g$ extends uniquely to an endomorphism of the $\N$-graded vertex algebra $V_{\Lie}$.
Furthermore, if $g$ is an automorphism, then
the extension is an automorphism.
\el

\begin{proof} Since $A\oplus B$ generates $V_{\Lie}$ as a vertex algebra,
the uniqueness is clear. 
It remains to prove the existence. By Lemma \ref{lend-lie}, 
we have a grading-preserving endomorphism $\hat{g}$ of 
the Lie algebra $\Lie$,
hence a grading-preserving endomorphism of the universal enveloping algebra $U(\Lie)$.
Consequently, $\hat{g}$ preserves the Lie subalgebra $\Lie^{<0}$ and 
its universal enveloping algebra $U(\Lie^{<0})$.
It follows from the construction of $V_{\Lie}$ that 
there exists a linear endomorphism $\bar{g}$ of
$V_{\Lie}$ such that $\bar{g}({\bf 1})={\bf 1}$ and 
$$\bar{g}(u_nv)=g(u)_n\bar{g}(v)$$
for $u\in A\oplus B$, $v\in V_{\Lie}$, $n\in \Z$. 
Since $V_{\Lie}$ is generated by $A\oplus B$, 
it follows (cf. [LLi]) that $\bar{g}$ is an endomorphism of $V_{\Lie}$. 
It is clear that $\bar{g}$ extends $g$.

If $g$ is an automorphism of the $1$-truncated conformal algebra $C=A\oplus B$,
{}from the first assertion we have vertex algebra endomorphisms $\bar{g}$ 
and $\overline{g^{-1}}$ of $V_{\Lie}$, extending $g$ and $g^{-1}$, respectively.
Since $gg^{-1}=g^{-1}g=1$ on $A\oplus B$ and since 
$A\oplus B$ generates $V_{\Lie}$ as a vertex algebra,
we have $\bar{g}\overline{g^{-1}}=\overline{g^{-1}}\bar{g}=1$.
Thus, $\bar{g}$ is an automorphism of $V_{\Lie}$.
\end{proof}

\bp{pem-B} 
Let $g$ be an endomorphism of a vertex $A$-algebroid $B$. 
Then $g$ extends uniquely to an endomorphism of $V_B$
as an $\N$-graded vertex algebra. Furthermore, if $g$ is an automorphism, then
the extension is an automorphism.
\ep

\begin{proof} 
The uniqueness is clear, as $A\oplus B$ generates $V_{B}$ as a vertex algebra.
For the the existence, first by Lemma \ref{lva}, 
we have a grading-preserving endomorphism $\bar{g}$ of 
the vertex algebra $V_{\Lie}$, extending $g$.
Now we show that $\bar{g}$ reduces to an endomorphism of $V_B$. 
Recall that $V_{B}=V_{\Lie}/I_{B}$, where $I_{B}$ is the two-sided ideal of
$V_{\Lie}$, generated by
$$E={\rm span}\{  e-{\bf 1},\; a(-1)a'-aa',\; a(-1)b-ab\;|\; a,a'\in A,\; b\in B\}.$$
Now, we must prove $\bar{g}(I_B)\subset I_B$. As $E$ generates $I_{B}$ as a two-sided ideal,
it suffices to prove that $\bar{g}(E)\subset E$.
Let $a,a'\in A,\; b\in B$. We have
\begin{eqnarray*}
& &\bar{g}(e-{\bf 1})=e-{\bf 1}\in E,\\
&&\bar{g}(a(-1)a'-aa')=g(a)(-1)g(a')-g(a)g(a')\in E,\\
&&\bar{g}(a(-1)b-ab)=g(a)(-1)g(b)-g(a)g(b)\in E.
\end{eqnarray*}
This proves $\bar{g}(E)\subset E$.
Therefore,  $\bar{g}$
reduces to an endomorphism of the $\N$-graded vertex algebra  $V_B$. 
The second assertion follows immediately from the proof
of the second assertion of Lemma \ref{lva}.
\end{proof}


Recall that $\Aut^{o} (V_B)$ denotes the group of grading-preserving automorphisms of $V_B$,
namely the group of automorphisms of $V_B$ as an $\N$-graded vertex algebra.
Combining Lemma \ref{lav} with Proposition \ref{pem-B}, we have:

\bt{taut-group} Let $A$ be a unital commutative associative algebra and 
let $B$ be a vertex $A$-algebroid.
The group $\Aut^{o}(V_B)$ of (grading-preserving) automorphisms of
the $\N$-graded vertex algebra $V_{B}$ is isomorphic to 
the group of automorphisms of vertex $A$-algebroid $B$ with the restriction map
as an isomorphism.
\et

\section{Classification of graded simple twisted $V_B$-modules}

In this section we construct and classify graded simple twisted
$V_B$-modules by exploiting a twisted analogue of the Lie algebra $\Lie$.
First, we recall the definition of the notion of twisted
module for a vertex algebra and we discuss several properties of
twisted modules.

Let $V$ be a vertex algebra and let $g$ be an automorphism of $V$ of 
order $T<\infty$. Decompose $V$ into eigenspaces of $g$:
$$V=\coprod_{r=0}^{T-1} V^r,
\ \ \text{ where }V^r=\{ v\in V\;  | \; g(v)=e^{2r\pi \sqrt{-1} /T}v\}.$$

A {\em $g$-twisted $V$-module} (see \cite{le}, \cite{flm},
\cite{ffr}, \cite{dong}) is a vector space $M$ equipped with a linear map 
\begin{eqnarray}
Y_M: V&\rightarrow& (\End M)[[x^{\frac{1}{T}},x^{-\frac{1}{T}}]]\nonumber\\
u&\mapsto& Y_{M}(u,x)=\sum_{n\in \frac{1}{T}\Z}u_nx^{-n-1}
\end{eqnarray}
satisfying the following conditions: 
\begin{enumerate}
\item For $u\in V,\; w\in M$, $u_n w=0$ for $n\in \frac{1}{T}\Z$ sufficiently large. 
\item $Y_M({\bf 1},x)=1_{M}$ (the identity operator on $M$). 
\item For $u\in V^r$ with $0\leq r\leq T-1$,
\begin{eqnarray}\label{etwisted-vform}
Y_M(u,x)=\sum_{n\in \frac{r}{T}+\Z}u_nx^{-n-1}
\in x^{-\frac{r}{T}}(\End M)[[x,x^{-1}]].
\end{eqnarray}
\item For $u\in V^r$ with $0\leq r\leq T-1$, $v\in V$,
\begin{eqnarray}\label{tji}
& &x_0^{-1}\delta\left(\frac{x_1-x_2}{x_0}\right) Y_M(u,x_1)Y_M(v,x_2)
-x_0^{-1}\delta \left(\frac{x_2-x_1}{-x_0}\right)Y_M(v,x_2)Y_M(u,x_1)\nonumber\\
& &
=x_2^{-1}\left(\frac{x_1-x_0}{x_2}\right)^{-r/T}\delta\left(\frac{x_1-x_0}{x_2}\right)
Y_M(Y(u,x_0)v,x_2)
\end{eqnarray}
(the {\em twisted Jacobi identity}).
\end{enumerate}

\br{rsubalgebra}
{\em  Let $(M,Y_{M})$ be a $g$-twisted $V$-module and 
let $U$ be any vertex subalgebra of $V^{0}$. Then $M$ is a $U$-module.
In particular, if $g$ is taken to be the identity map, the notion of $g$-twisted $V$-module reduces
to that of $V$-module while the twisted Jacobi identity reduces to
the ordinary (untwisted) Jacobi identity.}
\er

The following was proved in \cite{dlm-twisted} (cf. \cite{dlm-reg}):

\bl{ldlm-twisted}
Let $(M,Y_{M})$ be a $g$-twisted $V$-module. Then
\begin{eqnarray}
Y_M(\D v,x)=\frac{d}{dx}Y_M(v,x)
\end{eqnarray}
for $v\in V$, where $\D v=v_{-2}{\bf 1}$.
\el

\begin{rem} 
{\em For $v\in V$, $u\in V^r$, $p\in \Z$ and $s,t\in \Q$, 
comparing the coefficients of $z_0^{-p-1}z_1^{-s-1}z_2^{-t-1}$ on
the both sides of the twisted Jacobi identity (\ref{tji}) we get
\begin{eqnarray}\label{tji2}
\sum_{m\geq 0}{s\choose m}(u_{p+m}v)_{s+t-m}=\sum_{m\geq
0}(-1)^m{p\choose m}\{u_{p+s-m}v_{t+m}-(-1)^pv_{p+t-m}u_{s+m}\}.
\end{eqnarray}

By taking $\Res_{x_0}$ of (\ref{tji}), we obtain the {\em twisted
commutator formulae}:
\begin{eqnarray}\label{etwisted-comm}
& & [Y_M(u,x_1),Y_M(v,x_2)]\nonumber\\
&
&=\Res_{x_0}x_2^{-1}\left(\frac{x_1-x_0}{x_2}\right)^{-r/T}\delta\left(\frac{x_1-x_0}{x_2}\right)Y_M(Y(u,x_0)v,x_2).
\end{eqnarray}
Multiplying (\ref{tji}) by
$\left(\frac{x_1-x_0}{x_2}\right)^{\frac{r}{T}}$ and then taking
$\Res_{x_1}$, we obtain the {\em twisted iterate formulae}:
\begin{equation}\label{tif}
Y_M(Y(u,x_0)v,x_2)=\Res_{x_1}\left(\frac{x_1-x_0}{x_2}\right)^{\frac{r}{T}}
\cdot X
\end{equation}
where
$$X=x_0^{-1}\delta\left(\frac{x_1-x_2}{x_0}\right)Y_M(u,x_1)Y_M(v,x_2)
-x_0^{-1}\delta\left(\frac{x_2-x_1}{-x_0}\right)Y_M(v,x_2)Y_{M}(u,x_1).
$$
{}From the twisted Jacobi identity one has the following {\em twisted
weak associativity}: For $u\in V^r$ with $0\le r\le T-1$ and for $v\in V,\; w\in W$,
$$(x_0+x_2)^{k+\frac{r}{T}}Y_M(u,x_0+x_2)Y_M(v,x_2)w
=(x_2+x_0)^{k+\frac{r}{T}}Y_M(Y(u,x_0)v,x_2)w$$
where $k$ is a nonnegative integer such that
$x^{k+\frac{r}{T}}Y_M(u,x)w\in M[[x]]$.
One can prove (cf. \cite{li-twisted}; Lemma 2.8) that the twisted Jacobi identity 
is equivalent to
the twisted commutator formulae and the twisted weak associativity.}
\end{rem}

Let $M$ be a $g$-twisted $V$-module. For a subset $U$ of $M$, 
denote the smallest $g$-twisted $V$-submodule containing $U$ by
$\<U\>$, which is called the $g$-twisted $V$-submodule
generated by $U$. Just as with untwisted modules, 
{}from the twisted weak associativity, we have
$$\<U\>={\rm span}\{ v_{n}w \;| \;  v\in V,\; n\in \frac{1}{T}\Z,\;  w\in U\}.$$ 
We define
\begin{equation}
\Ann_V(U)=\{v\in V\; | \; Y(v,x)w=0\ \ \text{ for }w\in U\},
\end{equation}
the {\em annihilator of $U$ in $V$}.

\begin{prop}\label{iann} For any subset $U$ of a $g$-twisted $V$-module $M$, the
annihilator $\Ann_V(U)$ is an ideal of $V$. Moreover,
$$\Ann_V(U)=\Ann_V(\< U\>).$$
\end{prop}

\begin{proof} It follows immediately from the proof of Proposition 4.5.11 in [LLi]
with the weak associativity and the weak commutativity
relations being replaced  by the
twisted associativity and the twisted commutativity relations,
respectively. 
\end{proof}

Let $S$ be a subset of $V$. Define 
$$\Ann_M(S)=\{w\in M\;  |\; Y_M(v,x)w=0 \ \ \text{ for } v\in S\},$$ 
the {\em annihilator
of $S$ in $M$}. By suitably modifying the proof of
Proposition 4.5.14 in [LLi] and replacing the weak commutativity and
Proposition 4.5.11 (of [LLi]) in the proof of Proposition 4.5.14
in [LLi] by the twisted commutativity relation and Proposition
\ref{iann}, respectively, we immediately have:

\begin{prop}\label{mann} 
For a subset $S$ of $V$, the annihilator $\Ann_M(S)$
is a $g$-twisted $V$-submodule of $M$. Furthermore,
$$\Ann_M(S)=\Ann_M(\< S\>).$$
Here, $\<S\>$ is the ideal of $V$ generated by $S$.
\end{prop}

We shall use the following result of \cite{li-twisted} (Lemma 2.11):

\bl{lLi3}
Let $V$ be a vertex algebra with an automorphism $g$ of order $T$ and 
let $a\in V^{k},\; b,u^{0},\dots,u^{r}\in V$ with $0\le k\le T-1$.
If
\begin{eqnarray}
[Y(a,x_{1}),Y(b,x_{2})]=\sum_{j=0}^{r}\frac{1}{j!}Y(u^{j},x_{2})
\left(\frac{\partial}{\partial x_{2}}\right)^{j}x_{1}^{-1}\delta\left(\frac{x_{2}}{x_{1}}\right)
\end{eqnarray}
acting on $V$, then for any $g$-twisted $V$-module
$(M,Y_{M})$ we have
\begin{eqnarray}
[Y_{M}(a,x_{1}),Y_{M}(b,x_{2})]=\sum_{j=0}^{r}\frac{1}{j!}Y_{M}(u^{j},x_{2})
\left(\frac{\partial}{\partial x_{2}}\right)^{j}x_{1}^{-1}\delta\left(\frac{x_{2}}{x_{1}}\right)
\left(\frac{x_{2}}{x_{1}}\right)^{\frac{k}{T}}
\end{eqnarray}
acting on $M$. On the other hand, the converse is also true
for any faithful $g$-twisted $V$-module $(M,Y_{M})$. 
\el

\bd{dgraded-twisted-module}
{\em Let $V=\coprod_{m\in\Z}V_{(m)}$ be a $\Z$-graded vertex algebra.
A {\em
$\frac{1}{T}\N$-graded $g$-twisted $V$-module} is a $g$-twisted
$V$-module $M$ equipped with a $\frac{1}{T}\N$-grading
$$M=\coprod_{n\in\frac{1}{T}\N}M(n)$$
such that
 $$v_mM(n)\subset M(n+p-m-1)$$ 
for $v\in V_{(p)},\; m,n\in \frac{1}{T}\Z$ with $p\in \Z$.}
\ed

Next, we study $\frac{1}{T}\N$-graded $g$-twisted modules 
for the $\N$-graded vertex algebra $V_{B}$
associated to a vertex $A$-algebroid $B$, where $g$ is an automorphism of
order $T<\infty$ of the $\N$-graded vertex algebra $V_{B}$ and of the
vertex $A$-algebroid $B$ (cf. Theorem \ref{taut-group}).

Noticing that $A\oplus B$ is a $1$-truncated conformal algebra, 
we start with a general $1$-truncated conformal algebra $C=C_{0}\oplus C_{1}$
with an automorphism $g$ of $C$ of order $T<\infty$.
Associated with the $1$-truncated conformal algebra $C=C_{0}\oplus C_{1}$ 
we have the Lie algebra $\Lie(C)$ and the vertex algebra $V_{\Lie(C)}$ with
$C$ as a generating subspace. In view of Lemma \ref{lva} 
$g$ is an order-$T$ automorphism of the vertex algebra $V_{\Lie(C)}$.

The following result can be found in \cite{dlm-twisted}
(cf. \cite{bor}):

\bl{ldlm2}
Let $V$ be a vertex algebra and let $T$ be a positive integer.
Set
\begin{eqnarray}
L_{T}(V)=V\otimes \C[t^{\frac{1}{T}},t^{-\frac{1}{T}}],
\end{eqnarray}
a vector space, and set
$$\hat{\partial}=\D\otimes 1+1\otimes \frac{d}{dt},$$
a linear operator on $L_{T}(V)$. Then the bilinear (multiplicative) 
operation on $L_{T}(V)$, defined by
\begin{eqnarray}\label{etwisted-comm-lemma}
[u\otimes t^{m},v\otimes t^{n}]
=\sum_{i\ge 0}{m\choose
i}(u_{i}v\otimes t^{m+n-i})
\end{eqnarray}
for $u,v\in V,\; m,n\in \frac{1}{T}\Z$,
gives rise to a Lie algebra structure on $L_{T}(V)/\hat{\partial} L_{T}(V)$, which is denoted by
$\Lie(V,T)$. Furthermore, any order-$T$ automorphism $g$ of $C$ 
gives rise to an order-$T$ automorphism, also  denoted by $g$,
of $\Lie(V,T)$, where
\begin{eqnarray}
g(v\otimes t^{n})=e^{-2n\pi \sqrt{-1}}(gv\otimes t^{n})
\end{eqnarray}
for $v\in V,\; n\in \frac{1}{T}\Z$.
\el

Specializing Lemma \ref{ldlm2} with $V=V_{\Lie(C)}$, we have 
a Lie algebra $\Lie(V_{\Lie(C)},T)$ and an automorphism $g$.
For $u\in C,\; m\in \frac{1}{T}\Z$, denote by $u(m)$ the canonical
image of $u\otimes t^{m}$ in $\Lie(V_{\Lie(C)},T)$. We have
\begin{eqnarray}
& &(\partial a)(m)=-ma(m-1),\\
& &[u(m),v(n)]=\sum_{i=0}^{1}{m\choose i}(u_{i}v)(m+n-i)
\end{eqnarray}
for $a\in C_{0},\; u,v\in C,\; m,n\in \frac{1}{T}\Z$.
Because $u_{i}v\in C$ for $u,v\in C,\; i\ge 0$,
we see that $u(m)\; u\in C,\; m\in \frac{1}{T}\Z$ span
a Lie subalgebra $\Lie(C,T)$ of $\Lie(V_{\Lie(C)},T)$.
Denote by $\Lie(C,g)$ the $g$-fixed point Lie subalgebra:
\begin{eqnarray}
\Lie(C,g)=\Lie(C,T)^{g}.
\end{eqnarray}

Using Lemma \ref{ldlm2} we immediately have:

\begin{prop}
Let $C=C_{0}\oplus C_{1}$ be a $1$-truncated conformal algebra and
let $g$ be an order-$T$ automorphism of $C$. 
Then 
\begin{eqnarray}
\Lie(C,g)=L(C,g)/\hat \partial L(C_{0},g),
\end{eqnarray}
as a vector space, where
\begin{eqnarray}
L(C,g)=\coprod_{r=0}^{T-1}C^r \otimes t^{r/T}\C[t,t^{-1}],
\end{eqnarray}
$L(C_{0},g)$ is a subspace defined in the obvious way, and
$$\hat \partial=\partial\otimes 1+1\otimes d/dt: L(C_{0},g)\rightarrow L(C,g).$$
For $u\in C^{r}$ with $0\le r\le T-1$ and for $n\in \Z$, 
denote by $u(n+r/T)$
the canonical image of $u\otimes t^{n+r/T}$ in $\Lie(C,g)$.
Then the following relations hold
for $a\in C_{0}^r, \;a'\in C_{0}^{r'}, \;b\in C_{1}^s, \;b'\in C_{1}^{s'}, \; 
m,n\in \Z$:
\begin{eqnarray}
& &(\partial a)(m+r/T)=-(m+r/T)a(m-1+r/T),\\
& &[ a(m+r/T), a'(n+r'/T)]= 0, \label{e31}\\
& &{[ a(m+r/T), b(n+s/T)]} = (a_0b)(m+n+(r+s)/T), \label{e32}\\
& &{[ b(m+s/T), b'(n+s'/T)]}=(b_0b')(m+n+(s+s')/T)\nonumber\\
& &\hspace{3cm} \ \ \ \ \ \ +(m+s/T)(b_1b')(m+n+(s+s')/T-1).\label{e34}
\end{eqnarray}
\end{prop}


We define
\begin{eqnarray*}
& &\deg a(n+r/T)=-n-1\ \ \text{for }a\in C_{0}^{r}, \;n\in \Z,\\
& &\deg b(n+r/T)=-n\ \ \text{for }b\in C_{1}^{r},\; n\in \Z,
\end{eqnarray*}
making $\Lie(C,g)$ a $\frac{1}{T}\Z$-graded Lie algebra.
For $n\in \frac{1}{T}\Z$, denote by $\Lie(C,g)_{(n)}$ the degree-$n$
subspace. We have the following triangular
decomposition
$$\Lie(C,g)=\Lie(C,g)_+\oplus \Lie(C,g)_{(0)}\oplus \Lie(C,g)_{-},$$
where $\Lie(C,g)_{\pm}=\coprod_{0<n\in \frac{1}{T}\Z}\Lie(C,g)_{(\pm
n)}.$ Notice that $\Lie(C,g)_{(0)}$ is spanned by the elements
$a(-1),b(0)$ for $a\in C_{0}^0,\;b\in C_{1}^0$.

For $u\in C^{r}$ with $0\le r\le T-1$, form the generating function
\begin{eqnarray}
u(x)=\sum_{n\in \frac{r}{T}+\Z}u(n)x^{-n-1}\in
{\mathcal{L}}(C,g)[[x^{\frac{1}{T}}, x^{-\frac{1}{T}}]].
\end{eqnarray}
For any $\Lie (C,g)$-module $W$, we consider
$u(x)$ as an element of $(\End W)[[x^{\frac{1}{T}},x^{-\frac{1}{T}}]]$, 
which we denote by $u_{W}(x)$:
\begin{eqnarray}
u_W(x)=u(x)=\sum_{n\in \frac{1}{T}\Z}u(n)x^{-n-1}
\in (\End W)[[x^{\frac{1}{T}},x^{-\frac{1}{T}}]].
\end{eqnarray}

\bl{ltwisted-Lie-rep}
The commutation relations (\ref{e31})--(\ref{e34}) amount to 
the following relations in terms of generating functions:
\begin{eqnarray}
[a(x_1),a'(x_2)]&=&0,\label{etwisted-1}\\
{[a(x_1),b'(x_2)]}
&=&x_2^{-1}\left(\frac{x_1}{x_2}\right)^{-\frac{r}{T}}\delta\left(\frac{x_1}{x_2}\right)(a_0b')(x_2),\\
{[b(x_1),b'(x_2)]}
&=&x_2^{-1}\left(\frac{x_1}{x_2}\right)^{-\frac{r}{T}}\delta\left(\frac{x_1}{x_2}\right)(b_0b')(x_2)
\nonumber\\
& &+
(b_1b')(x_2)\frac{\partial}{\partial{x_2}}x_2^{-1}\left(\frac{x_1}{x_2}\right)^{-\frac{r}{T}}
\delta\left(\frac{x_1}{x_2}\right)\label{etwisted-3}
\end{eqnarray}
for $a\in C_{0}^r,\;b\in C_{1}^r, \;a'\in C_{0}$, and $b'\in C_{1}$. 
Moreover, we have
\begin{eqnarray}
(\partial a)(x)=\frac{d}{dx}a(x)\ \ \ \mbox{ for }a\in C_{0}.
\end{eqnarray}
\el

{}From these relations we immediately have:

\bc{c33} For $a,a'\in C_{0},\;b,b'\in C_{1}$,
\begin{eqnarray}
[a(x_1),a'(x_2)]&=&0,\\
(x_1-x_2)[a(x_1),b(x_2)]&=&0,\\
(x_1-x_2)^2[b(x_1),b'(x_2)]&=&0.
\end{eqnarray}
\ec

\begin{de} {\em An $\mathcal{L}(C,g)$-module $W$ is said to be {\em restricted}
if for any $w\in W,\;u\in C^{r}$ with $0\le r\le T-1$,
$u(n+r/T)w=0$ for $n\in \Z$ sufficiently large,
that is, $u_{W}(x)\in \Hom (W,W((x^{\frac{1}{T}})))$ for $u\in C$.}
\end{de}

The following result is analogous to a result of 
\cite{li-twisted} for twisted affine Lie algebras:

\begin{prop}\label{ptwisted-affine} 
Let $C=C_{0}\oplus C_{1}$ be a $1$-truncated conformal algebra
and let $g$ be an automorphism of $V_{\Lie(C)}$ of order $T$, which is extended from
an automorphism of $C$.
Every $g$-twisted $V_{\Lie(C)}$-module $W$ is naturally a restricted
$\Lie(C,g)$-module with $u_W(x)=Y_W(u,x)$ for $u\in C$.
Moreover, the set of $g$-twisted $V_{\Lie(C)}$-submodules of $W$ is precisely the set of
$\Lie(C,g)$-submodules of $W$.
On the other hand, for any restricted $\Lie(C,g)$-module $W$, there exists a unique
$g$-twisted $V_{\Lie(C)}$-module structure $Y_{W}$ on $W$ such that
\begin{eqnarray}
Y_W(u,x)=u_W(x)\;\;\;\mbox{ for }u\in C=C_{0}\oplus C_{1}\subset V_{\Lie(C)}.
\end{eqnarray}
\end{prop}

\begin{proof} On the vertex algebra $V_{\Lie(C)}$, the following
relations hold for $a,a'\in C_{0},\;b,b'\in C_{1}$:
\begin{eqnarray}
& &[Y(a,x_1),Y(a',x_2)]=0,\\
& &[Y(a,x_1),Y(b',x_2)]
=x_2^{-1}\delta\left(\frac{x_1}{x_2}\right)Y(a_0b',x_2),\\
& &[Y(b,x_1),Y(b',x_2)]
=x_2^{-1}\delta\left(\frac{x_1}{x_2}\right)Y(b_0b',x_2)
+Y(b_1b',x_2)\frac{\partial}{\partial{x_2}}x_2^{-1}
\delta\left(\frac{x_1}{x_2}\right).\ \ \ \ \ \ \ \
\end{eqnarray}
{}From Lemmas \ref{lLi3} and \ref{ltwisted-Lie-rep}, 
every $g$-twisted $V_{\Lie(C)}$-module $W$ is naturally a restricted
$\Lie(C,g)$-module with $u_W(x)=Y_W(u,x)$ for $u\in C$.
As $C$ generates $V_{\Lie(C)}$ as a vertex algebra,
the set of $g$-twisted $V_{\Lie(C)}$-submodules of $W$ 
is precisely the set of $\Lie(C,g)$-submodules of $W$.

Let $S={\rm span}\{u_W(x)\;|\;u\in C\}$. 
In view of Corollary \ref{c33}, $S$ is a local subspace of 
$\Hom (W,W((x^{\frac{1}{T}})))$.
Note that 
$\Hom (W,W((x^{\frac{1}{T}})))$ is naturally $\Z/T\Z$-graded:
$$\Hom (W,W((x^{\frac{1}{T}})))
=\coprod_{r=0}^{T-1}x^{\frac{r}{T}}\Hom (W,W((x)))$$
and $S$ is a graded subspace. Let $\sigma_{T}$ be the linear
automorphism of
$\Hom (W,W((x^{\frac{1}{T}})))$ defined by
$$\sigma_{T}(\alpha(x))=e^{-2r\pi\sqrt{-1}/T}\alpha(x)$$
for $\alpha(x)\in x^{\frac{r}{T}}\Hom (W,W((x)))$ with $0\le r\le T-1$
(cf. (\ref{etwisted-vform})).

{}From \cite{li-twisted}, $S$ generates a vertex
algebra $\<S\>$ inside $\Hom (W,W((x^{\frac{1}{T}})))$ with the
identity operator $1_{W}$ as the vacuum vector and 
with $\sigma_{T}$ as an automorphism.
Furthermore, $W$ is naturally a faithful 
$\sigma_{T}$-twisted $\<S\>$-module with 
$Y_{W}(\alpha(x),x_{0})=\alpha(x_{0})$. 
With the relations (\ref{etwisted-1})-(\ref{etwisted-3}),
by Lemma \ref{lLi3}, we have
\begin{eqnarray}
& &[Y(a_{W}(x),x_1),Y(a'_{W}(x),x_2)]=0,\\
& &[Y(a_{W}(x),x_1),Y(b'_{W}(x),x_2)]
=x_2^{-1}\delta\left(\frac{x_1}{x_2}\right)Y((a_0b')_{W}(x),x_2),\\
& &[Y(b_{W}(x),x_1),Y(b'_{W}(x),x_2)]
=x_2^{-1}\delta\left(\frac{x_1}{x_2}\right)Y((b_0b')_{W}(x),x_2)\nonumber\\
& &\ \ \ \ \hspace{2cm} +Y((b_1b')_{W}(x),x_2)\frac{\partial}{\partial{x_2}}x_2^{-1}
\delta\left(\frac{x_1}{x_2}\right)
\end{eqnarray}
for $a\in C_{0}^r,\;b\in C_{1}^r, \;a'\in C_{0}$, and $b'\in C_{1}$.
We also have
$$Y((\partial a)_{W}(x),x_{1})=Y\left(\frac{d}{dx} a_{W}(x),x_{1}\right)
=\frac{\partial}{\partial x_{1}}
Y(a_{W}(x),x_{1})$$
for $a\in A$. By Lemmas \ref{lLi3} and \ref{ltwisted-Lie-rep}, 
$\<S\>$ is naturally an $\Lie(C)$-module with
$u_{\<S\>}(x_{1})=Y(u_{W}(x),x_{1})$ for $u\in C$.
Furthermore, $\<S\>$ as an $\Lie(C)$-module is generated by $1_{W}$
and we have $u_{W}(x)_{n}1_{W}=0$ for $u\in C,\; n\ge 0$.
{}From the construction of $V_{\Lie(C)}$ as an $\Lie(C)$-module,
there exists a unique $\Lie(C)$-module 
homomorphism $\psi$ from $V_{\Lie(C)}$ 
to $\<S\>$, sending ${\bf 1}$ to $1_{W}$. As $V_{\Lie(C)}$ as a vertex
algebra is generated by $C$, $\psi$ is a vertex algebra homomorphism.
We have
$$\psi(u)=\psi(u(-1){\bf 1})=u_{W}(x)_{-1}1_{W}=u_{W}(x)$$
for $u\in C$. It is clear that
$\psi(C^{r})\subset S^{r}$ for $0\le r\le T-1$. As $C$ generates
$V_{\Lie(C)}$ as a vertex algebra, $\psi$ preserves the
$\Z/T\Z$-gradings, i.e., $\sigma_{T}\psi =\psi g$.
Consequently, $W$ is a $g$-twisted $V_{\Lie(C)}$-module.
\end{proof}

For the rest of this paper, {\em we assume that 
$A$ is a unital commutative associative algebra
whose identity is denoted by $e$ and $B$ is a vertex $A$-algebroid and 
we assume that $g\in \Aut^{o} (V_B)$ with $o(g)=T<\infty$.} 
Recall that $C=A\oplus B$ is naturally a $1$-truncated conformal algebra.
An $\Lie (C,g)$-module of {\em level} $k\in\C$ is an
$\Lie (C,g)$-module on which $e(-1)$ acts as scalar $k$.

Immediately from Proposition \ref{ptwisted-affine} we have:

\bp{psubmodules} 
Every $g$-twisted $V_B$-module is naturally a restricted
$\Lie(C,g)$-module of level $1$.
Moreover, the set of $g$-twisted $V_B$-submodules is
precisely the set of $\Lie(C,g)$-submodules.
\ep

We have the following decompositions into $g$-eigenspaces:
$$V_{B}=\coprod_{r=0}^{T-1}V_{B}^r$$ 
and
\begin{eqnarray*}
A=\coprod_{r=0}^{T-1}A^r\ \ \text{ and}\ \ \ \ 
B=\coprod_{r=0}^{T-1}B^r.
\end{eqnarray*}
Clearly, $A^0$ is a subalgebra of $A$, containing the identity, 
and $B^0$ is a
vertex $A^0$-algebroid. Furthermore, by Lemma \ref{br2},
$B^0/A^0\partial A^0$ is a Lie $A^0$-algebroid.
Set
\begin{eqnarray}
I=\sum_{r=1}^{T-1}A^{r}\cdot A^{T-r}\subset A^{0}.
\end{eqnarray}
It is clear that $I$ is a two-sided ideal of $A^{0}$, so that 
$A^{0}/I$ is a unital commutative associative algebra. Furthermore, 
$B^{0}/(I\cdot B^{0}+A^0\partial A^0)$ is a Lie $A^0/I$-algebroid.

\begin{prop} Let $M=\coprod_{n\in \frac{1}{T}\N}M(n)$ 
be a $\frac{1}{T}\N$-graded $g$-twisted $V_B$-module.
Then $M(0)$ is a module for the Lie $A^0$-algebroid $B^0/A^0\partial A^0$
with
\begin{eqnarray}
& & a\cdot w=a_{-1}w\ \ \ \ \mbox{ for }a\in A^0,\; w\in M(0),\\
& &b\cdot w=b_0w\ \ \ \ \mbox{ for }b\in B^0, \; w\in M(0).
\end{eqnarray}
Furthermore, for $a\in A^r$,  $a'\in A^{T-r}$, $b\in B^{T-r}$ with $0<r\leq T-1$
and for $w\in M(0)$, we have 
$(aa')\cdot w=0$ and $(ab)\cdot w=(1-\frac{r}{T})(a_0b)\cdot w$.
\end{prop}

\begin{proof} Let $U$ be the vertex subalgebra of $V_{B}$ 
generated by $A^{0}\oplus B^{0}$.
As $A^{0}\oplus B^{0}\subset V_{B}^{0}$, 
$U$ is actually a vertex subalgebra of $V_{B}^{0}$.
{}From Remark \ref{rsubalgebra}, $M$ is a $U$-module. 
With $(V_{B})_{(0)}=A$ and $(V_{B})_{(1)}=B$, we have
$(V_{B}^{0})_{(0)}=A^{0}$ and $(V_{B}^{0})_{(1)}=B^{0}$.
Consequently, we have $U_{(0)}=A^{0}$ and $U_{(1)}=B^{0}$.
It follows from the construction of $V_{B^{0}}$ that $U$ is a 
homomorphic image of the vertex algebra $V_{B^{0}}$, so that
$W$ is naturally a $V_{B^{0}}$-module.
By \cite{liy} (Proposition 4.8), 
$W(0)$ is naturally a module for the Lie $A^0$-algebroid $B^0/A^0\partial A^0$.

Let $a\in A^r,\; a'\in A^{T-r},\;b\in B^{T-r},\; w\in M(0)$ with
$0<r\le T-1$.
By substituting $u=a$, $v=a'$, $p=-1$, $s=-1+\frac{r}{T}$,
$t=-\frac{r}{T}$ in (\ref{tji2}), we get
\begin{eqnarray*}
(aa')\cdot w&=&(aa')_{-1}w\\
&=&(a(-1)a')_{-1}w\\
&=&\sum_{m\geq 0}a_{-2+\frac{r}{T}-m}a'_{-\frac{r}{T}+m}w+a'_{-1-\frac{r}{T}-m}a_{-1+\frac{r}{T}+m}w\\
&=&0.
\end{eqnarray*}
Similarly, by substituting $u=a$, $v=b$, $p=-1$,
$s=-1+\frac{r}{T}$ and $t=1-\frac{r}{T}$ in (\ref{tji2}), we have
\begin{eqnarray*}
(ab)\cdot w&=&(ab)_0w\\
&=&(a(-1)b)_0w\\
&=&(1-\frac{r}{T})(a_0b)_{-1}w
+\sum_{m\geq 0}\{a_{-2+\frac{r}{T}-m}b_{1-\frac{r}{T}+m}+b_{-\frac{r}{T}-m}a_{-1+\frac{r}{T}+m}\}w\\
&=&(1-\frac{r}{T})(a_0b)\cdot w,
\end{eqnarray*}
completing the proof.
\end{proof}

Let $U$ be a module for the Lie $A^0$-algebroid $B^0/A^0\partial A^0$ such
that $(aa')\cdot u=0$ and $(ab)\cdot u=(1-\frac{r}{T})(a_0b)\cdot
u$ for $a\in A^r,\; a'\in A^{T-r}, \;b\in B^{T-r},\;u\in U, \;0<r\leq
T-1$. We are going to construct a $\frac{1}{T}\N$-graded
$g$-twisted $V_B$-module $M=\coprod_{n\in \frac{n}{T}\N}M(n)$ with
$M(0)=U$ as a module for the Lie $A^0$-algebroid
$B^0/A^0\partial A^0$.

First, $U$ is a module for the Lie algebra $A^0\oplus B^0/\partial
A^0$. Recall that $\Lie(C,g)_{(0)}=A^0\oplus B^0/\partial A^0$. For
convenience, we set $\Lie (C,g)_{\leq 0}=\Lie (C,g)_{(0)}\oplus \Lie
(g)_{-}$. Then $U$ is an $\Lie(C,g)_{\leq 0}$-module under the
following actions
\begin{eqnarray*}
a(n+\frac{r}{T}-1)\cdot u&=&\delta_{n+\frac{r}{T},0}au,\\
b(n+\frac{r}{T})\cdot u&=& \delta_{n+\frac{r}{T},0}bu
\end{eqnarray*}
for $a\in A^r,\;b\in B^r,\;n\geq 0$. Next, we form the induced
$\Lie(C,g)$-module
\begin{eqnarray}
M_g(U)=\Ind_{\Lie(C,g)_{\leq 0}}^{\Lie(C,g)}U
=U(\Lie(C,g))\otimes_{U(\Lie(C,g)_{\leq 0})}U.
\end{eqnarray}
We endow $U$ with degree $0$, making $M_g(U)$ a
$\frac{1}{T}\N$-graded restricted $\Lie(C,g)$-module.
By Proposition \ref{ptwisted-affine},
$M_g(U)$ is naturally a $g$-twisted $V_{\Lie(C)}$-module.
In view of the P-B-W theorem, we may and we should consider 
$U$ as the degree-zero subspace of $M_g(U)$.

We set
\begin{eqnarray}
W_g(U)={\rm span}\{ v_nu \; | \; v\in E, \; n\in \frac{1}{T}\Z,\;  u\in U\}
\subset M_g(U)
\end{eqnarray}
and define
\begin{eqnarray}
M_{B}(U)=M_g(U)/U(\Lie(C,g))W_g(U).
\end{eqnarray}
Since $U(\Lie(C,g))W_g(U)$ 
is an $\Lie(C,g)$-submodule of $M_g(U)$,
by Proposition \ref{ptwisted-affine} $U(\Lie(C,g))W_g(U)$ is a $g$-twisted
$V_{\Lie(C)}$-submodule. Then $M_B(U)$ is a $g$-twisted
$V_{\Lie(C)}$-module. Clearly, $M_B(U)$ is generated by $\bar{U}$
the image of $U$ in $M_{B}(U)$. In fact,
$M_{B}(U)$ is a $g$-twisted $V_{B}$-module by the following:

\bl{genu} Let $(M, Y_M)$ be a $g$-twisted $V_{\Lie}$-module.
Suppose that for $a\in
A^r,\; a'\in A,\; b\in B$ with $0\leq r\leq T-1$,
\begin{eqnarray}
Y_M(e,x)w&=&w,\label{ee1}\\
Y_M(a(-1)a',x)w&=&Y_M(aa',x)w,\label{ea(-1)a'}\\
Y_M(a(-1)b,x)w&=&Y_M(ab,x)w\label{ea(-1)b}
\end{eqnarray}
for all $w\in K$, where $K$ is a generating subspace of $M$.
Then $M$ is naturally a $g$-twisted $V_B$-module.
\el

\begin{proof} Recall that $$E={\rm span}\{e-{\bf 1}, \; a(-1)a'-aa',\;
a(-1)b-ab \; | \; a,a'\in A,b\in B\}\subset V_{\Lie(C)}.$$ By
(\ref{ee1})-(\ref{ea(-1)b}), we have $K\subset \Ann_M(E)$. By using
Propositions \ref{mann}, we have $\Ann_M(I_B)=\Ann_M(E)$. Since
$\Ann_M(I_B)$ is a $g$-twisted $V_{\Lie(C)}$ submodule of $M$ and $M$
is generated by $K$, we have $\Ann_M(I_B)=M$. This implies that $M$
is a $g$-twisted $V_B$-module.
\end{proof}

One can see that Lemma \ref{genu} indeed implies that
$M_B(U)$ is naturally a $g$-twisted $V_B$-module.
Furthermore we have:

\bt{tuniversal} Let $U$ be a module for the Lie $A^0$-algebroid $B^0/A^0\partial
A^0$ such that 
\begin{eqnarray}
(aa')\cdot u=0\ \ \mbox{ and }(ab)\cdot
u=(1-\frac{r}{T})(a_0b)\cdot u
\end{eqnarray}
for $a\in A^r,\; a'\in A^{T-r},\; b\in
B^{T-r},\; u\in U,\; r\neq 0$. Then $M_B(U)$ is naturally a $g$-twisted
$V_B$-module such that $M_B(U)(0)=U$.
\et

\begin{proof} To show that $M_B(U)(0)=U$,
we must prove that $(U(\Lie(g))W_g(U))(0)=0$.
First we show that $W_g(U)(0)=0$. 
Notice that for $v\in (V_{\Lie(C)})_{(m)}^{r}$ with $m\in \Z$,
we have $\deg (v_{k+r/T})=m-k-r/T-1$ for $k\in \Z$.
Then from the definition of $W_g(U)$,
$W_g(U)(0)$ is spanned by the vectors
$$(e-{\bf 1})_{-1}u,\;\; (a(-1)a')_{-1}u-(aa')_{-1}u,\;\;
(a(-1)b)_0u-(ab)_0u$$ 
for $u\in U$, $a\in A^{r},\;a'\in A^{T-r}$, $b\in B^{T-r}$ with $0\le
r\le T-1$. Since $e_{-1}$ acts as $e$ (the identity of $A^{0}$) on $U$,
we have $(e-{\bf 1})_{-1}u=0$ for $u\in U$.
If $r=0$, by (\ref{tji2}), we have
$$(a_{-1}a')_{-1}u=a(-1)a'(-1)u=a(a'u)=(aa')u=(aa')_{-1}u, $$ and
$$(a(-1)b)_0u=a(-1)b(0)u=a(bu)=(ab)u=(ab)_0u.$$
Next, we assume that $r>0$. By (\ref{tji2}), we have
\begin{eqnarray*}
(a(-1)a')_{-1}u&=&\sum_{i=0}^{\infty}a(-1-i+\frac{r}{T})a'(i-1-\frac{r}{T})u
+\sum_{i=0}^{\infty}a'(-2-i-\frac{r}{T})a(i+\frac{r}{T})u\\
&=&a(-1+\frac{r}{T})a'(-1-\frac{r}{T})u\\
&=&a'(-1-\frac{r}{T})a(-1-\frac{r}{T})u\\
&=&0\\
&=&(aa')_{-1}u
\end{eqnarray*}
and
\begin{eqnarray*}
& &(a(-1)b)_{0}u\\
&=&\sum_{i=0}^{\infty}a(-1-i+\frac{r}{T})b(i-\frac{r}{T})u
+\sum_{i=0}^{\infty}b(-i-1-\frac{r}{T})a(i+\frac{r}{T})u
-\frac{r}{T}(a_0b)_{-1}u\\
&=&a(-1+\frac{r}{T})b(-\frac{r}{T})u-\frac{r}{T}(a_0b)\cdot u\\
&=&b(-\frac{r}{T})a(-1+\frac{r}{T})u+(a_0b)_{-1}u-\frac{r}{T}(a_0b)_{-1}u\\
&=&(1-\frac{r}{T})(a_0b)\cdot u\\
&=&(ab)\cdot u.
\end{eqnarray*}
Hence, $W_g(U)(0)=0$.

Next, we show that $\Lie (C,g)_{\leq 0}W_g(U)\subset W_g(U)$. 
Recall from \cite{liy} (Lemma 4.2) that
$$v_{i}E\subset E\;\;\;\mbox{ for }v\in C=A\oplus B,\;i\ge 0. $$
As $M_{g}(U)$ is a $\frac{1}{T}\N$-graded $\Lie (C,g)$-module with $U$
as the degree-zero subspace, we have that $\Lie (C,g)_{\leq 0}U\subset U$.
For $v\in C=A\oplus B,\; c\in E,\; m,t\in \frac{1}{T}\Z,\;u\in
U$, {}from the twisted commutator formula (\ref{etwisted-comm}) 
(cf. (\ref{etwisted-comm-lemma})), we have
$$v_{m}c_tu=c_tv_{m}u+\sum_{i\geq 0}{m\choose i}(v_ic)_{m-t-i}u.$$
These immediately imply
that $\Lie(C,g)_{\leq 0}W_g(U)\subset W_g(U)$.
Then
\begin{eqnarray*}
U(\Lie(C,g))W_g(U)&=&U(\Lie(C,g)_{+})U(\Lie(C,g)_{\leq
0})W_g(U)\\
&=&U(\Lie(C,g)_{+})W_g(U)\\
&=&W_g(U)+\Lie(C,g)_{+}U(\Lie(C,g)_{+})W_g(U),
\end{eqnarray*}
which implies that $(U(\Lie(C,g))W_g(U))(0)=0$.
This completes the proof.
\end{proof}

Next, we continue with Theorem \ref{tuniversal}
to construct and classify $\frac{1}{T}\N$-graded simple 
$g$-twisted $V_B$-modules.
Let $U$ be a module for the Lie $A^0$-algebroid $B^0/A^0\partial
A^0$ as in Theorem \ref{tuniversal}. 
Let $J(U)$ be the sum of all graded $\Lie(C,g)$-submodules
of $M_g(U)$ with trivial degree-zero subspaces.
Then $J(U)$ is the unique maximal graded $\Lie(C,g)$-submodule
of $M_g(U)$ with the property that $J(U)\cap U=0$.
Set
\begin{eqnarray}
L_g(U)=M_g(U)/J(U),
\end{eqnarray}
a $\frac{1}{T}\N$-graded $g$-twisted $V_{B}$-module.

\bl{lsimple-module} Let $U$ be a module for the Lie $A^0$-algebroid
$B^0/A^0\partial A^0$ as in Theorem \ref{tuniversal}.
Then $L_g(U)$ is a
$\frac{1}{T}\N$-graded $g$-twisted $V_B$-module such that
$L_g(U)(0)=U$ as a module for the Lie $A^0$-algebroid
$B^0/A^0\partial A^0$ and such that for any nonzero graded
submodule $W$ of $L_{g}(U)$, we have $W(0)\ne 0$.
Furthermore, if $U$ is a simple
$B^0/A^0\partial A^0$-module, $L_g(U)$ is a 
$\frac{1}{T}\N$-graded simple $g$-twisted $V_B$-module.
\el

\begin{proof} It is similar to the proof of Theorem 4.12 in \cite{liy}.
\end{proof}

\begin{lem} Let $W=\coprod_{n\in\frac{1}{T}\Z}W(n)$ be a
$\frac{1}{T}\N$-graded simple $g$-twisted $V_B$-module with $W(0)\neq 0$. Then
$W\cong L_{g}(W(0))$.
\end{lem}

\begin{proof} It is similar to the proof of Lemma 4.13 in
\cite{liy}.
\end{proof}

To summarize we have:

\begin{thm} Let $H$ be a complete set of 
equivalence class representatives of simple modules for the Lie $A^0$-algebroid
$B^0/A^0\partial A^0$ satisfying the condition that
$$ (aa')U=0,\; \; ((ab)-(1-\frac{r}{T})(a_0b))U=0$$
for $a\in A^r, a'\in A^{T-r}, b\in B^{T-r}$ with $0<r<T$.
Then $\{
L_g(U)\;|\;U\in H\}$ is a complete set of 
equivalence class representatives of $\frac{1}{T}\N$-graded simple
$g$-twisted $V_B$-modules.
\end{thm}

\begin{proof} It is similar to the proof of Theorem 4.14 in
\cite{liy}.
\end{proof}

\br{rimplication} 
{\em Taking $g=1$ the identity map of $V_{B}$, we recover
Theorem 4.14 of \cite{liy}:
If $H$ is a complete set of 
equivalence class representatives of simple modules for the Lie $A$-algebroid
$B/A\partial A$, then $\{L_{1}(U)\;|\;U\in H\}$ is a complete
set of equivalence class representatives of 
$\N$-graded simple $V_B$-modules.}
\er


\begin{thebibliography}{BFFLS}

\bibitem [BBS]{bbs}
S. Berman, Y. Billig and J. Szmigielski,
Vertex operator algebras and the representation theory of toroidal algebras,
in: {\em Recent
Developments in Infinite-Dimensional Lie Algebras and Conformal Field
Theory (Charlottesville, VA, 2000)}, Contemporary Math. {\bf 297},
Amer. Math. Soc., Providence, 2002, 1-26.

\bibitem [BDT]{bdt}
S. Berman, C. Dong and S. Tan, Representations of a class of lattice
type vertex algebras, {\em J. Pure Applied Algebra}
{\bf 176} (2002), 27-47.

\bibitem [B]{bor}
R. E. Borcherds, Vertex algebras, Kac-Moody algebras, and the Monster,
{\it Proc. Natl. Acad. Sci. USA} {\bf 83} (1986), 3068-3071.

\bibitem[Br1]{br1}
P. Bressler, Vertex algebroids I, arXiv: math.AG/0202185.

\bibitem[Br2]{br2}
P. Bressler, Vertex algebroids II, arXiv: math.AG/0304115.

\bibitem[D]{dong}
C. Dong, Twisted modules for vertex algebras associated with
even lattices, {\em J. Algebra} {\bf 165} (1994), 91-112.

\bibitem [DLM1]{dlm-reg}
C. Dong, H.-S. Li and G. Mason, 
Regularity of rational vertex operator algebras,
{\em Advances in Math.} {\bf 132} (1997), 148-166.

\bibitem[DLM2]{dlm-twisted}
C. Dong, H.-S. Li and G. Mason, Twisted representations of vertex operator 
algebras, {\em Math. Annalen} {\bf 310} (1998), 571-700.

\bibitem[DLM3]{dlm-vpa}
C. Dong, H.-S. Li and G. Mason, Vertex Lie
algebra, vertex Poisson algebras and vertex algebras, in: {\em Recent
Developments in Infinite-Dimensional Lie Algebras and Conformal Field
Theory}, Proceedings of an International Conference at University of
Virginia, May 2000, Contemporary Math. {\bf 297} (2002), 69-96.

\bibitem[FFR]{ffr}
A. J. Feingold, I. B. Frenkel and J. F. X. Ries, {\it Spinor
Construction of Vertex Operator Algebras, Triality, and
$E_{8}^{(1)}$}, Contemporary Math. {\bf 121}, 1991.

\bibitem[FB]{fbbook}
E. Frenkel and D. Ben-Zvi, {\em Vertex Algebras and Algebraic Curves},
Mathematical Surveys and Monographs, Vol. 88, Amer. Math. Soc.,
Providence, 2001.

\bibitem[FF1]{ff88}
B. Feigin and E. Frenkel, A family of representations of affine Lie
algebras, {\em Russian Math. Surveys} {\bf 43} (1988), 221-222.

\bibitem[FF2]{ffflag}
B. Feigin and E. Frenkel, Affine Kac-Moody algebras and semi-infinite
flag manifolds, {\em Commun. Math. Phys.} {\bf 128} (1990), 161-189.

\bibitem[FF3]{ff-representation}
B. Feigin and E. Frenkel, Representations of affine Kac-Moody Lie
algebras and bosonization, in: {\em Physics and Mathematics of
Strings}, World Scientific, 1990, 271-316.

\bibitem[FLM]{flm}
I. B. Frenkel, J. Lepowsky and A. Meurman, {\it Vertex Operator
Algebras and the Monster,} Pure and Applied Math., Vol. 134,
Academic Press, Boston, 1988.

\bibitem[GMS]{gms} V. Gorbounov, F. Malikov and V. Schechtman,
Gerbes of chiral differential operators, II, Vertex algebroids,
{\em Invent. Math}. {\bf 155} (2004), 605-680.

\bibitem[Le]{le}
J. Lepowsky, Calculus of twisted vertex operators, {\em Proc. Natl.
Acad. Sci. USA} {\bf 82} (1985), 8295-8299.

\bibitem[LLi]{lli}
J. Lepowsky and H.-S. Li, {\em Introduction to Vertex Operator Algebras
and Their Representations,} Progress in Math. {\bf 227},
Birkh$\ddot{a}$user, Boston, 2003.

\bibitem[Li1]{li1}
H.-S. Li, Representation theory and tensor product theory for vertex
operator algebras, Ph.D. thesis, Rutgers University, 1994.

\bibitem[Li2]{li-local}
H.-S. Li, Local systems of vertex operators, vertex superalgebras
and modules,
{\em J. Pure Appl. Alg.} {\bf 109} (1996), 143-195;  hep-th/9406185.

\bibitem[Li3]{li-twisted}
H.-S. Li, Local systems of twisted vertex operators, vertex
operator superalgebras and twisted modules, in: {\em Moonshine,
the Monster and Related Topics, Proc. Joint Summer Research
Conference, Mount Holyoke}, 1994, ed. by C. Dong and G. Mason,
Contemporary Math. {\bf 193}, Amer. Math. Soc., Providence, 1996,
203-236.

\bibitem[LY]{liy}
H.-S. Li and G. Yamskulna, On certain vertex algebras and their
modules associated with vertex algebroids, 
{\em J. Algebra} {\bf 283} (2005), 367-398.

\bibitem[MS]{ms}
F. Malikov and V. Schechtman, Chiral de Rham complex, II,
axXiv: math.AG/9901065.

\bibitem[MSV]{msv}
F. Malikov, V. Schechtman and A. Vaintrob, Chiral de Rham complex,
axXiv: math.AG/9803041.

\bibitem [W]{wakimoto}
M. Wakimoto, Fock representations of affine Lie algebra $A_{1}^{(1)}$,
{\em Commun. Math. Phys.} {\bf 104} (1986), 605-609.

\end{thebibliography}
\end{document}